\documentclass{amsart}

\usepackage{amsmath,amsthm,mathrsfs,amsrefs}
\usepackage{amssymb,amsfonts,tikz-cd}
\usepackage{url,relsize,xcolor}
\usepackage[hidelinks]{hyperref}
\usepackage{bbm}

\newcommand{\bbB}{\mathbb B}
\newcommand{\bbC}{\mathbb C}
\newcommand{\bbP}{\mathbb P}
\newcommand{\bbQ}{\mathbb Q}
\newcommand{\ZFC}{\mathsf{ZFC}}
\DeclareMathOperator{\Fix}{Fix}
\DeclareMathOperator{\cof}{cof}
\DeclareMathOperator{\Gen}{Gen}
\DeclareMathOperator{\RO}{RO}
\DeclareMathOperator{\Aut}{Aut}
\DeclareMathOperator{\Add}{Add}

\newtheorem{lemma}{Lemma}
\newtheorem{proposition}{Proposition}
\newtheorem{corollary}{Corollary}

\theoremstyle{definition}
\newtheorem{definition}{Definition}

\newtheorem*{ack}{Acknowledgments}
\newtheorem{example}{Example}

\theoremstyle{remark}
\newtheorem{remark}{Remark}

\author{Ben De Bondt}
\date{}
\title{A family of forcings with hyperfinite equivalence of generics}

\begin{document}
	
\maketitle	

\begin{abstract}
We identify a class of forcings that includes several classical forcings and prove  that for all forcings $\bbP$ in this class the associated \emph{equivalence of generics relation}  $\equiv^\bbP_M$, as defined by Smythe in \cite{Smythe}, is a hyperfinite CBER.
\end{abstract}	

For $M$ a countable transitive model of $\ZFC$ \footnote{Or a sufficiently large fragment thereof.}
(ctm) and $\bbP$ a forcing notion in~$M$, define the following Borel equivalence relation $\equiv^\bbP_M$ on the standard Borel space $\Gen_{M}^\bbP$ of all filters that are $\bbP$-generic over $M$:
	\[ G \equiv^\bbP_M H  \Leftrightarrow M[G] = M[H].\]
Equivalently (up to Borel isomorphism), by the Vopenka-Hajek theorem, $\equiv^\bbP_M$ can also be thought of as the orbit equivalence relation induced by the Borel action $\Aut^M(\bbB) \times \Gen_{M}^\bbB \to \Gen_{M}^\bbB: (\alpha,G) \to \alpha[G]$ of the countable group $\Aut^M(\bbB)$ on the standard Borel space $\Gen_{M}^\bbB$, where $\bbB$ is the Boolean completion of $\bbP$ in~$M$ (\cite[Lemma~2.7, Theorem~2.16]{Smythe}).
These countable Borel equivalence relations (CBERs) were introduced and first systematically studied in \cite{Smythe}.
A central open problem in the theory of CBERs asks whether every hyper-hyperfinite CBER (i.e.\ countable increasing union of hyperfinite CBERs) is again hyperfinite. Since for Cohen forcing $\bbC = \Add(\omega,1)^M$ the CBER  $\equiv^\bbC_M$ is hyper-hyperfinite and is not known to be hyperfinite, it is natural to ask whether a counterexample to this increasing unions problem can be found amongst the CBERs $\equiv^\bbP_M$. In particular, Smythe asked (\cite[Question 5.1.]{Smythe}):
\\
\emph{what can we say about the complexity of $\equiv^\bbP_M$ when $\bbP$ is,
e.g., Sacks, Miller, \hbox{Mathias}, or Laver forcing, or any of the other classical
forcing notions? In particular, is $\equiv^\bbP_M$ hyperfinite in any of these cases?}
\\
For the specific case of Prikry forcing, hyperfiniteness was established in \cite{Calderoni_Sinapova}. We identify here a large class of forcing notions for which $\equiv^\bbP_M$ is hyperfinite. This answers Smythe's question for Sacks\footnote{Hyperfiniteness in the case of Sacks forcing was known prior to this note (\cite{Sinapova}).}, Miller and Laver forcing and several other classical forcing notions such as Namba forcing and all forcings that don't add reals. Our methods do not (directly) apply to the case $\bbP = \Add(\omega,1)^M$.

\begin{definition}\label{def1}\,
	\begin{enumerate}
		\item Let $\bbP$ be a forcing notion and $\bbQ$ a complete subforcing of $\bbP$. We call $\bbQ$ \emph{absorbing} in $\bbP$ if $1_\bbP \Vdash_{\bbP} \dot{G} \in V[\dot{G} \cap \bbQ].$
		\item 	Let $\mathcal{C}$ be a definable class of forcings (e.g.\ the class of ccc forcings). We say that a complete Boolean algebra $\bbB$ has the \emph{$\mathcal{C}$-absorption property} if \[ 1_\bbB \Vdash_{\bbB} \text{ ``if $\bbB / (\dot{G} \cap \bbC) \in \mathcal{C}^{V[\dot{G} \cap \bbC]},$ then $\dot{G} \in V[\dot{G} \cap \bbC]$'',}\] for every complete subalgebra $\bbC$ of $\bbB$.
	\end{enumerate}
\end{definition}

So, re-expressing from the point of view of the complete subforcing, $\bbQ$ being absorbing in $\bbP$ simply means that 
\[
1_\bbQ\Vdash_{\bbQ}
\text{``$\bbP/\dot G_{\bbQ}$ is forcing trivial,''}
\]
(where we call a forcing notion \emph{forcing trivial} if it never adds a new set) and $\bbB$ has the $\mathcal{C}$-absorption property if for all of its complete subalgebras $\bbC$:
	\[
	1_\bbC\Vdash_{\bbC} \bbB/\dot G_{\bbC} \in \dot{\mathcal{C}} \Rightarrow
	\text{$\bbB/\dot G_{\bbC} $ is forcing trivial.}
	\]
	
\begin{example}
(The Boolean completions of) Sacks, Miller and Laver forcing all have the ccc-absorption property.	More generally, any complete Boolean algebra $\bbB$ which is nowhere ccc and is such that any $\bbB$-forcing extension is a \emph{minimal} extension has the ccc-absorption property. Indeed, consider $V\subseteq V[\mathbb{C}\cap G]\subseteq V[G]$ with $G$ generic for $\mathbb{B}$ over $V$ and $\mathbb{C}$ a complete subalgebra of $\mathbb{B}$ and assume that $\mathbb{B}/(\mathbb{C}\cap G)$ is ccc in $V[G\cap\mathbb{C}]$. Suppose towards a contradiction that  $V[G\cap \mathbb{C}] \neq V[G]$.   By minimality, then $V=V[G\cap \mathbb{C}]$, hence $V[G]$ would be realised as a generic extension of $V[G\cap \mathbb{C}]$ both by $\mathbb{B}$ and a ccc forcing, contradicting that $\bbB$ is nowhere ccc.
\end{example}

\begin{example}
Similarly, also Miller-style Namba forcing has the ccc-absorption property, again the argument exploits a minimality property of the respective forcing. 
Let $\kappa\geq\omega_2$ be regular, let $\mathbb{B}=\RO(\mathbb{N}_\kappa)$. Consider $V\subseteq V[\mathbb{C}\cap G]\subseteq V[G]$ with $G$ generic for $\mathbb{B}$ over $V$ and $\mathbb{C}$ a complete subalgebra of~$\mathbb{B}$. Assume that $\mathbb{B}/(\mathbb{C}\cap G)$ is ccc in $V[G\cap\mathbb{C}]$. Since ccc forcing does not change cofinalities we must have that $\cof^{V[G\cap\mathbb{C}]}(\kappa)=\omega$, but then the ``$\cof(\kappa)=\omega$''-minimality of $\mathbb{N}_\kappa$ (\cite[Lemma 2.2]{Levine}) gives that indeed $V[G\cap\mathbb{C}]=V[G]$.
\end{example}
	
When $\Gamma$ is a set of automorphisms of a complete Boolean algebra $\bbB$, we write $\Fix(\bbB|\Gamma)$ for the set of elements of $\bbB$ that are fixed by all elements of $\Gamma$. Note that then $\Fix(\bbB|\Gamma)$ is a complete subalgebra of $\bbB$.

\begin{definition}
For $\bbB$ a complete Boolean algebra and $\Gamma \leq \Aut(\bbB)$, we define (by transferred epithet) the group $\Gamma$ to be an \emph{absorbing} subgroup of $\Aut(\bbB)$ if $\Fix(\bbB | \Gamma)$ is absorbing in $\bbB$. We define $\bbB$ to be \emph{spongy} if every finitely generated subgroup $\Gamma$ of $\Aut(\bbB)$ is absorbing.
\end{definition}

Note that if $\bbC$ is a complete subalgebra of $\bbB$ that is absorbing in $\bbB$, then every complete subalgebra $\mathbb{D}$ of $\bbB$ containing $\bbC$ is again absorbing in $\bbB$. Hence if $\Gamma \leq \Aut(\bbB)$ is absorbing, then also every subgroup of $\Gamma$ is absorbing.

\begin{lemma}
Let $\bbB$ be a complete Boolean algebra in the ctm $M$. The following two are equivalent:
\begin{enumerate}
	\item $M \vDash \text{``$\bbB$ is spongy''}$,
	\item  there exists (in $V$) an increasing sequence $(\Gamma_n:n<\omega)$ of absorbing subgroups of $\Aut^M(\bbB)$ such that $\Aut^M(\bbB) = \bigcup_n \Gamma_n$ and for all $n$, $\Gamma_n \in M$  and $M\vDash  \text{``$\Gamma_n$ is absorbing''}$.
\end{enumerate}
\end{lemma}

\begin{proof}
	$(1) \Rightarrow (2)$ follows because $M$ is countable, hence $\Aut^M(\bbB)$ can be written as a countable increasing union of its finitely generated subgroups.
Conversely, $(2)$ implies that in $M$ every finitely generated subgroup $\Gamma$ of $\Aut^M(\bbB)$ is contained in a subgroup of $\Aut^M(\bbB)$ that is absorbing, and therefore that $\Gamma$ is absorbing as well.
\end{proof}

 The following is our central lemma (and motivates the foregoing peculiar definitions).

\begin{lemma}\label{lem1}
	Let $\bbB$ be a complete Boolean algebra in the ctm $M$, and let $(\Gamma_n:n<\omega) \in V$ be an increasing sequence of subgroups of $\Aut^M(\bbB)$ such that $\Aut^M(\bbB) = \bigcup_n \Gamma_n$ and for every $n< \omega$:
	
	\begin{itemize}
		\item $\Gamma_n \in M$,
		\item $M\vDash  \text{``$\Gamma_n$ is absorbing''}$,
	\end{itemize}
	then the Borel equivalence relation $\equiv_M^\bbB$ is hyperfinite.\\	
	In particular, if, in $M$, the Boolean algebra $\bbB$ is spongy, then  $\equiv_M^\bbB$ is hyperfinite.
\end{lemma}
 \begin{proof}
	Following the set-up of \cite{Smythe}, let $\Gen^\bbB_M$ be the standard Borel space of $\mathbb{B}$-filters generic over $M$. 
	Consider for each $n$ the equivalence relation $E_n$ on $\Gen_M^\bbB$ defined by:
	\[ G\, E_n\, H :\iff G\cap\Fix(\bbB|\Gamma_n)=H\cap \Fix(\bbB|\Gamma_n).\] 
	Note that:
	\begin{itemize}
		\item $E_n\subseteq E_{n+1}\subseteq\, \equiv_M^\bbB$ for every $n$, because $\Fix(\bbB|\Gamma_n)$ is absorbing in $\bbB$ so $G\cap \Fix(\bbB|\Gamma_n)=H\cap\Fix(\bbB|\Gamma_n)$ implies
		\[ M[G]=M[G\cap \Fix(\bbB|\Gamma_n)]=M[H\cap\Fix(\bbB|\Gamma_n)]=M[H],   \]
		\item the equivalence relations $E_n$ are smooth, because $G\mapsto G \cap\Fix(\bbB|\Gamma_n)$ gives a reduction of $E_n$ to equality in the standard Borel space $2^\bbB$. 
	\end{itemize}
	Moreover, the relations $E_n$ for $n<\omega$ exhaust $\equiv_M^\bbB$. To see this, recall that by the Vop\v{e}nka-H\'ajek theorem (see \cite[Theorem~2.16]{Smythe}) we have that $\equiv_M^\bbB$ coincides with the orbit equivalence relation induced by the action of $\Aut^M(\bbB)$ on $\Gen_M^\bbB$ given by $H^\gamma:=\gamma[H]$ for $H\in\Gen_M^\bbB$, $\gamma\in\Aut^M(\bbB)$.
	Hence, if $G\equiv_M^\bbB H$ then $G=\gamma[H]$ for some $\gamma\in\Aut^M(\bbB)$, which gives 
	\[ G\cap \Fix(\bbB|\Gamma_n) = H\cap \Fix(\bbB|\Gamma_n)  \]
	for all $n<\omega$ so large that $\gamma\in\Gamma_n$.
	
	It follows that $\equiv_M^\mathbb{B}$ is an increasing union of smooth countable Borel equivalence relations and by \cite[Theorem~5.1]{Dougherty_Jackson_Kechris} it follows that $\equiv_M^\mathbb{B}$ is hyperfinite.
\end{proof}

We next show that the Boolean completions of all of the following forcings are spongy and therefore for each of these forcings the associated generics-equivalence-relation is hyperfinite:

\begin{itemize}
	\item every forcing which does not add any reals,
	\item (Miller-style) Namba forcing,
	\item Sacks forcing,
	\item Laver forcing,
	\item Miller forcing.
\end{itemize}

We first handle forcings that don't add reals. Note that this includes many important examples of forcings, including all countably closed forcings and other $\omega$-Baire forcings (e.g. shooting a club through a stationary set using countable conditions, Souslin trees, $\ldots$) and various types of Prikry forcings (we thus generalize a result from \cite{Calderoni_Sinapova}).

\begin{lemma}\label{lem3}
	Let, in $M$, $\bbB$ be a complete Boolean algebra and $\Gamma$ a subgroup of $\Aut^M(\bbB)$. Then \[1_\bbB \Vdash_{\bbB} \bbB / (\dot{G} \cap \Fix(\bbB | \Gamma)) = \bigcup_{\gamma\in\Gamma}\gamma[\dot{G}].\]
	In particular, if $\Gamma$ is countable, then
	$1_\bbB \Vdash_{\bbB} ``\bbB / (\dot{G} \cap \Fix(\bbB | \Gamma)) \text{ is $\sigma$\text{-centered}}$''.
\end{lemma}
\begin{proof}
	Let $G$ be generic for $\bbB$ over $M$ and let $\bbC=\Fix(\bbB|\Gamma)$. Recall that $\bbB/(G \cap \bbC) = \{b \in \bbB \colon(\forall c \in G \cap \bbC)\ b \wedge c \neq 0\}$. 
	It is clear that  $\bbB / (G \cap \bbC) \supseteq \bigcup_{\gamma\in\Gamma}\gamma[G]$. For the other direction, let $b \in \bbB$ be compatible with all elements of $G \cap \bbC$. Then $\bigvee_{\gamma\in\Gamma} \gamma(b) \in \bbC$ is compatible with all elements of $G \cap \bbC$, hence $\bigvee_{\gamma\in\Gamma} \gamma(b)\in G\cap \bbC$ and therefore $b\in \bigcup_{\gamma\in\Gamma} \gamma[G]$. (Alternatively, note that $\bigvee_{\gamma\in\Gamma}\gamma(b)=\bigwedge\{ c\in\bbC\colon b\leq c  \}=\pi(b)$ with $\pi:\mathbb{B}\to\bbC$ the induced projection map.)
\end{proof}

\begin{lemma}\label{lem2}
	Let, in $M$, $\bbC$ be a complete subalgebra of a complete Boolean algebra~$\bbB$. Let $G$ be $\bbB$-generic over $M$. If  $M[G] \vDash ``\bbB / (G \cap \bbC) \text{ is ccc}$'', then 
	\[ (2^\omega)^{M[G\cap \bbC]} = (2^\omega)^{M[G]} \Leftrightarrow M[G\cap \bbC] = M[G].\]
\end{lemma}

\begin{proof}Suppose that $M[G\cap\mathbb{C}]$ and $M[G]$ have the same reals.  Let $\mathbb{D}$ in $M[G\cap\mathbb{C}]$ be the Boolean completion of $\bbB / (G \cap \bbC)$. Because $M[G] \vDash ``\mathbb{D} \text{ is ccc}$'' and $\omega_1^{M[G]}=\omega_1^{M[G\cap\mathbb{C}]}$, we can choose $d \in \mathbb{D}$ such that $M[G]$ is a $\mathbb{D} \upharpoonright d$-extension and either (a) $\mathbb{D} \upharpoonright d$ is atomic or (b) $\mathbb{D} \upharpoonright d$ is atomless and has ccc square. Then, since they also have the same reals, $M[G]$ is a forcing extension of $M[G\cap \bbC]$ by forcing with a complete Boolean algebra which is either (a) trivial or (b) is both atomless, $\omega$-Baire and has ccc square.  Since the square of an atomless Souslin algebra is never ccc, the latter is impossible and we must have $ M[G\cap \bbC] = M[G]$. 
\end{proof}

\begin{proposition}
	If $M\models``\bbB\text{ adds no reals}$'', then $\bbB$ is spongy in $M$ and by consequence $\equiv_M^\bbB$ is hyperfinite.
\end{proposition}

\begin{proof}
In $M$, let $\Gamma$ be any countable subgroup of $\Aut^M(\bbB)$, we show $\Gamma$ is absorbing. Let $G$ be generic for $\bbB$ over $M$. Let $\bbC= \Fix(\bbB | \Gamma)$. By Lemma~\ref{lem3}, $M[G] \vDash ``\bbB / ({G} \cap \Fix(\bbB | \Gamma))$ is $\sigma$-centered'', hence by Lemma~\ref{lem2}, $ M[G\cap \bbC] = M[G]$. This shows that $\Fix(\bbB|\Gamma)$ is absorbing in~$\bbB$. Hyperfiniteness of  $\equiv_M^\bbB$ now follows from Lemma~\ref{lem1}.
\end{proof}

Note that if $M$ satisfies $\mathsf{CH}$, then Namba forcing over $M$ doesn't add reals, so by the foregoing proposition in this case  $\equiv_M^\bbB$ is hyperfinite (where $\bbB$ is the Boolean completion of Namba forcing). We will show now that for Miller-style Namba forcing Lemma~\ref{lem1} can be applied to get hyperfiniteness without any additional assumptions on~$M$. The argument actually applies to all $\omega_1$-preserving forcing notions with the ccc-absorption property (in the sense of Definition~\ref{def1}(2)).

\begin{proposition}\label{Prop }
	If the complete Boolean algebra $\bbB$ has the ccc-absorption property in $M$ and preserves $\omega_1$, then $\bbB$ is spongy in $M$ and $\equiv^\bbB_M$ is a hyperfinite Borel equivalence relation.
\end{proposition}

\begin{corollary}
(The Boolean completions of) Miller-style Namba forcing, Sacks forcing, Laver forcing and Miller forcing are all spongy and therefore the corresponding Borel equivalence relations $\equiv^\bbB_M$ are hyperfinite.
\end{corollary}
 
\begin{proof}[Proof of Proposition~\ref{Prop }]
	 Suppose $\mathbb{B}$ has the ccc-absorption property in $M$ and doesn't collapse $\omega_1$. We show that, in $M$, for every countable subgroup $\Gamma$ of $\Aut^M(\mathbb{B})$, $\Fix(\mathbb{B}|\Gamma)$ is absorbing in $\mathbb{B}$.
	 Let $\Gamma\in M$ be a countable subgroup of $\Aut^M(\mathbb{B})$. Define $\mathbb{C}=\Fix(\mathbb{B}|\Gamma)$ and consider $M\subseteq M[\mathbb{C}\cap G]\subseteq M[G]$ with $G$ generic for $\mathbb{B}$ over $M$. Since by Lemma~\ref{lem3} the quotient forcing $\mathbb{B}/(\mathbb{C}\cap G)$ is ccc in $M[G]$ and $\omega_1^{M[G\cap\mathbb{C}]}=\omega_1^{M[G]}$, $\mathbb{B}/(\mathbb{C}\cap G)$ is also ccc in $M[G\cap\mathbb{C}]$. By the ccc-absorption property of $\mathbb{B}$, it follows that $M[G\cap \mathbb{C}]=M[G]$.
 \end{proof} 
 
 \begin{remark}
 Clearly, (the Boolean completion of) Cohen forcing $\Add(\omega,1)$ is \emph{not} spongy. Indeed, every minimal homeomorphism of Cantor space induces an automorphism of its regular open algebra that fixes only $0$ and $1$.
 \end{remark}
 
 \begin{ack}
 	This research was funded by the Deutsche Forschungsgemeinschaft (DFG, German
 	Research Foundation) under Germany’s Excellence Strategy EXC 2044/2-390685587,
 	Mathematics Münster: Dynamics–Geometry–Structure.
 	
This note grew out of an earlier proof attempt by the author to prove hyperfiniteness of $\equiv^\bbP_M$ for Namba forcing (inspired by the Prikry case in \cite{Calderoni_Sinapova} and work by Andreas Lietz on Mathias characterisations for Namba generics), which got strongly simplified through discussions with ChatGPT~5.6. In particular, the system pointed out in this discussion the elementary identity $\bbB/(G \cap \Fix(\bbB | \Gamma)) = \bigcup_{\gamma\in\Gamma}\gamma[G]$ which underlies the above arguments. The note is entirely human-written, the author is responsible for correctness of the content.
 \end{ack}

\end{document}